\documentclass{elsarticle}

\usepackage{geometry}

\usepackage{amssymb,amsmath,amsthm}
\usepackage{epsfig,graphicx,epstopdf}

\newtheorem{lem}{Lemma}[section]
\newtheorem{theorem}{Theorem}[section]
\newtheorem{remark}{Remark}[section]
\theoremstyle{example}
\newtheorem{ex}{Example}[section]

\allowdisplaybreaks
\begin{document}

\begin{frontmatter}

\title{Legendre wavelet collocation method combined with the Gauss--Jacobi quadrature for solving fractional delay-type integro-differential equations}

\author{S. Nemati$^{a,}$\footnote{Corresponding author. Email address: s.nemati@umz.ac.ir}, P.M. Lima$^b$, S. Sedaghat$^c$}
\address{$^a$Department of Mathematics, Faculty of Mathematical Sciences, University of Mazandaran, Babolsar, Iran}
\address{$^b$Centro de Matem\'{a}tica Computacional e Estoc\'astica, Instituto Superior T\'{e}cnico, Universidade de Lisboa, Av. Rovisco Pais, 1049-001 Lisboa, Portugal}
\address{$^c$Department of Mathematics, Buein Zahra Technical University, Buein Zahra, Qazvin, Iran}

\begin{abstract}
In this work, we present a collocation method based on the Legendre wavelet combined with the Gauss--Jacobi quadrature formula for solving a class of fractional delay-type integro-differential equations. The problem is considered with either initial or boundary conditions and the fractional derivative is described in the Caputo sense.  First, an approximation of the unknown solution is considered in terms of the Legendre wavelet basis functions. Then, we substitute this approximation and its derivatives into the considered equation. The Caputo derivative of the unknown function is approximated using the Gauss--Jacobi quadrature formula. By collocating the obtained residual at the well-known shifted Chebyshev points, we get a system of nonlinear algebraic equations. In order to obtain a continuous solution, some conditions are added to the resulting system. Some error bounds are given for the Legendre wavelet approximation of an arbitrary function. Finally, some examples are included to show the efficiency and accuracy of this new technique.

\end{abstract}

\begin{keyword}
Fractional delay-type integro-differential equations \sep Legendre wavelet \sep Gauss--Jacobi quadrature \sep Chebyshev collocation points

%% PACS codes here, in the form: \PACS code \sep code

\MSC 34A08 \sep 45J05 \sep 42C40
%% or \MSC[2008] code \sep code (2000 is the default)

\end{keyword}

\end{frontmatter}
%%%%%%%%%%%%%%%%%%%%%%%%%%%%%%%%%%%%%%%%%%%%%%%%%%%%%%%%%%%%%%%%%%%%%%%%%%%%%%%%%%%%%%%%%%%%%%%%%%%%%%%%%%%%%%
\section{Introduction}
\label{sec:1}
In recent decades, fractional calculus, which includes theories of derivatives and integrals of any arbitrary non-integer order, has been used for describing many real world phenomena, in such different domains as wave mechanics, viscoelasticity, physiology  or epidemiology  \cite{Das,Benson,Popovic,Sun,Babaei,Azadi}. Nevertheless, obtaining analytic solutions for fractional differential and integro-differential equations  is very difficult. Thus, many researchers have introduced and developed numerical methods in order to get numerical solutions for these two classes of equations (see for example  \cite{Arikoglu,Dehghan,Momani,Odibat,Saadatmandi,Nemati1,Nemati2}). Moreover, partial differential equations of fractional order have also been investigated and various numerical methods for their approximation have been introduced recently \cite{Badr,Singh}.

Delay differential and integro-differential equations are used in the modeling of a large number of problems in different fields of sciences such as ecology and epidemiology \cite{Baker}, immunology \cite{Fowler}, physiology \cite{Cooke} and  electrodynamics \cite{Dehghan1}. Fractional delay differential equations have been used to model certain processes and systems with memory and heritage properties. Discussions about the existence of the solution of such equations can be found in \cite{Lakshmikantham,Liao,Ye}. Because of the computational complexity of fractional derivatives in the delay case, the exact analytical solution of the fractional delay differential equations is hardly available. Therefore, in the last decades many researchers have been attracted to deal with the numerical solution of this class of problems (see for example \cite{Ordokhani,Rahimkhani,Saeed1,Yang,Morgado,Nemati3}).

Wavelets are a set of functions built by dilation and translation of a single function $\varphi(t)$ which is called the mother wavelet. These functions are known as a powerful tool to get the numerical solution of equations. The primary idea of wavelets (translation and  dilation) returns to the early 1960’s \cite{Chui,Daubechies}. Some developments exist concerning the multiresolution analysis algorithm based on wavelets \cite{Daubechies1} and the construction of compactly supported orthonormal wavelet bases \cite{Mallat}. Wavelets form an unconditional (Riesz) basis for $L^2(\mathbb{R})$ which is the space of all square integrable functions on the real line, i.e. any function in $L^2(\mathbb{R})$ can be decomposed and reconstructed in terms of wavelets \cite{Daubechies2}. Some researchers have constructed and used a number of wavelets for solving fractional differential equations such as B-spline wavelets \cite{Li1}, Haar wavelets \cite{Chen1}, Chebyshev wavelets \cite{Li2}, Legendre wavelets \cite{Jafari} and Bernoulli wavelets \cite{Ordokhani}. The advantage of wavelets in comparison with other basis functions is that when using them we can improve the accuracy of the approximation in two different ways: 1) by increasing the degree of the mother function (assuming that it is a polynomial); 2) by increasing the level of resolution, i.e. reducing the support of each basis function. Unfortunately no all the researchers who use wavelets for approximating fractional equations take this advantage; for example, in \cite{Saeed} the wavelets are used only with resolution level $k=1$, and only the degree of the polynomial is changed. 

In this paper, we consider the following fractional order nonlinear delay-type integro-differential equation
\begin{equation}\label{1.1}
^cD^\alpha y(x)=f\left(x,y(x),y'(x),y(px-\tau),y'(px-\tau),\int_{px-\tau}^xg(x,s,y(s))ds\right),
\end{equation}
subject to the initial conditions
\begin{equation}\label{1.2}
y(x)=y_0(x),\quad \quad x\leq 0, \quad y'(0)=y'_0,
\end{equation}
or boundary conditions
\begin{equation}\label{1.3}
y(x)=y_0(x),\quad x\leq 0,\quad y(l)=y_1,
\end{equation}
where $0\leq x \leq l$, $1<\alpha \leq 2$, $^cD^\alpha$ denotes the Caputo derivative, $f$ and $g$ are continuous linear or nonlinear functions and $y'_0$, $y_1$ and $p$ are constants. The parameter $\tau$ is called the delay component and is always non-negative. The delay component may be just constant (the constant delay case) or a function of $x$, $\tau=\tau(x)$ (the variable dependent delay case). The term $px-\tau$ is called delayed argument, $y(x)$ is the unknown function to be determined and $y(px -\tau)$ is the value of this function corresponding to the delayed argument. The function $y_0(x)$ is defined and continuously differentiable in $[\rho,0]$, where
$$\rho=\min_{x\geq 0}(px-\tau)=-\tau;\text{ moreover  $y_0'(0)=y'_0$}. $$
Delay problems of this type include as particular case the pantograph equation, when $p\ne 1$  and $\tau=0$. The existence and uniqueness of solution for this kind of problems has been investigated by many authors. In the non-delay case (when there is not delayed argument) we have just a nonlinear fractional equation; the existence and uniqueness of solution for such problems has been discussed for example in \cite{Diethelm}; among other restrictions the function $f$ is supposed to  satisfy the Lipschitz condition. In the case of  delayed argument, authors usually apply the method of steps to reduce the problem to a non-delay one; then existence and uniqueness is analyzed for the new problem. This was done, for example, in \cite{Morgado}, for the case of linear fractional equations with delayed arguments.

We present a collocation spectral method based on the Legendre wavelet basis functions combined with the Gauss--Jacobi quadrature formula to solve the above problem. We underline that the application of wavelets as basis functions has special advantages when problems of this class are considered. Actually in many cases the solution is not regular, in the sense that its first or higher derivatives may have discontinuities. In such cases, increasing the degree of polynomial approximation is not a solution to improve the accuracy of the approximation; however such improvement can be obtained by increasing the level of resolution. This is why in the numerical examples in Section \ref{sec:5} we vary both the polynomial degree and the resolution level of wavelets, and we check the effect on the quality of the approximation. To this aim, the rest of the paper is organized as follows: in Section \ref{sec:2}, we give some preliminaries of fractional calculus and introduce the Gauss--Jacobi quadrature formula and the Legendre wavelet basis functions. Section \ref{sec:3} is devoted to suggesting a numerical method based on the Legendre wavelet and using the Gauss--Jacobi quadrature for solving the fractional order delay-type integro-differential equations. In Section \ref{sec:4}, some error bounds are given. Section \ref{sec:5} includes some numerical examples. Finally, conclusions are given in Section \ref{sec:6}.
%%%%%%%%%%%%%%%%%%%%%%%%%%%%%%%%%%%%%%%%%%%%%%%%%%%%%%%%%%%%%%%%%%%%%%%%%%
\section{Basic concepts}\label{sec:2}
In this section, we present some preliminaries which will be used further in the next section.
\subsection{Preliminaries of fractional calculus}
Here, we briefly give some preliminaries and notations of fractional calculus. Two of the most important definitions for fractional integral and derivative operators are Riemann-Liouville integral and Caputo derivative. The Riemann-Liouville fractional integral operator $I^\alpha$ of order $\alpha\geq 0$ is defined as follows \cite{Podlubny}:
\begin{equation*}\label{0.1}
I^\alpha y(x)=\left\{\begin{array}{ll}
\frac{1}{\Gamma(\alpha)}\int_0^x (x-s)^{\alpha-1}y(s) ds,&\alpha>0,\\
y(x),& \alpha=0,
\end{array}
\right.
\end{equation*}
where $\Gamma(\alpha)$ is the gamma function defined by
\begin{equation*}\label{0.2}
\Gamma(\alpha)=\int_0^\infty x^{\alpha-1}e^{-x}dx.
\end{equation*}
Also, the Caputo fractional derivative operator $^cD^{\alpha}$ of order $\alpha$ is defined as follows  \cite{Podlubny}:
\begin{equation}\label{2.11}
^cD^{\alpha}y(x)=\frac{1}{\Gamma(n-\alpha)}\int_0^x (x-s)^{n-\alpha-1}y^{(n)}(s)ds,\quad n-1<\alpha\leq n,\quad n\in \mathbb{N},
\end{equation}
where $n=\lceil \alpha\rceil$ is the smallest integer greater than or equal to $\alpha$.
\subsection{Gauss--Jacobi quadrature formula}
The Jacobi polynomials $\{J_n^{(\lambda,\nu)}\}_{n=0}^{\infty}$, $\lambda,\nu>-1$, $x\in[-1,1]$, built a set of functions which are orthogonal with respect to the weight function 
\begin{equation*}
w^{\lambda,\nu}(x)=(1-x)^\lambda(1+x)^\nu.
\end{equation*}
The orthogonality property is satisfied for these polynomials as follows:
\begin{equation*}
\int_{-1}^{1}w^{\lambda,\nu}(x)J_n^{(\lambda,\nu)}(x)J_m^{(\lambda,\nu)}(x)dx=h_n^{\lambda,\nu}\delta_{mn},
\end{equation*}
where $\delta_{mn}$ is the Kronecker delta and 
\begin{equation*}
h_n^{\lambda,\nu}=\frac{2^{\lambda+\nu+1}\Gamma(\lambda+n+1)\Gamma(\nu+n+1)}{n!(\lambda+\nu+2n+1)\Gamma(\lambda+\nu+n+1)}.
\end{equation*}

The Gauss--Jacobi quadrature formula is given by
\begin{equation}\label{2.7}
\int_{-1}^{1}(1-x)^{\lambda}(1+x)^{\nu}f(x)dx=\sum_{s=1}^{N}w_sf(x_s)+R_N(f),
\end{equation}
where $x_s$, $s=1,\ldots,N$, are the roots of $J_N^{(\lambda,\nu)}$. Moreover, $w_s$, $s=1,\ldots,N$, are the weights and $R_N(f)$ is the remainder term which are, respectively, given by (see \cite{Shen}):
\begin{equation}\label{2.8}
w_s=\frac{2^{\lambda+\nu+1}\Gamma(\lambda+N+1)\Gamma(\nu+N+1)}{N!\Gamma(\lambda+\nu+N+1)(\frac{d}{dx}J_N^{(\lambda,\nu)}(x_s))^2(1-x_s^2)},
\end{equation}  
\begin{equation*}
\begin{split}
R_n(f)=& \frac{2^{\lambda+\nu+2N+1}N!\Gamma(\lambda+N+1)\Gamma(\nu+N+1)\Gamma(\lambda+\nu+N+1)}{(\lambda+\nu+2N+1)(\Gamma(\lambda+\nu+2N+1))^2}\\
&\quad\quad\quad\quad\quad\quad\quad\quad\quad\quad\quad\quad\quad\quad\times \frac{f^{(2N)}(\xi)}{(2N)!},\quad \xi\in(-1,1).
\end{split}
\end{equation*}
Taking into account the remainder term, formula \eqref{2.7} is exact for all polynomials of degree less than or equal to $2N-1$ and is valid if $f$ has no singularity in $(-1,1)$. We use this formula further in our numerical method.

\subsection{Legendre wavelets}
We define Legendre wavelet functions on the interval $[0,l)$ as
\begin{equation}\label{2.2}
\psi_{n,m}(x)=\left\{\begin{array}{ll}
2^{\frac{k}{2}}\sqrt{\frac{2m+1}{2l}}P_m\left(\frac{2^k}{l}x-2n+1\right),& \frac{n-1}{2^{k-1}}l\leq x<\frac{n}{2^{k-1}}l,\\
0,&\text{otherwise},
\end{array}
\right.
\end{equation}
where $k = 1, 2, 3, \ldots $, is the level of resolution, $n = 1, 2, 3,\ldots, 2^{k-1}$, $x$ is the normalized time, and $m = 0,1, 2, \cdots$, is the degree of the Legendre polynomial which is defined on the interval $[-1,1]$ and is given by the aid of the following recursive formula:
\begin{equation*}
P_m(x)=2P_{m-1}(x)-P_{m-2}(x),\quad m=2,3,\ldots,
\end{equation*}
with $P_0(x)=1$ and $P_1(x)=x$. Furthermore, in the special case of Jacobi polynomials with $\lambda=\nu=0$,  the Legendre polynomials are given. These polynomials are orthogonal with respect to the weight function $w(x)=1$, that is
\begin{equation*}
\int_{-1}^1P_i(x)P_j(x)dx=\frac{2}{2i+1}\delta_{ij}.
\end{equation*}

The first and second derivatives of the Legendre wavelet function $\psi_{n,m}(x)$ are given, respectively, by
\begin{equation}\label{2.9}
\psi_{n,m}^{'}(x)=\left\{\begin{array}{ll}
\frac{2^{\frac{3k}{2}}}{l}\sqrt{\frac{2m+1}{2l}}P_m^{'}\left(\frac{2^k}{l}x-2n+1\right),& \frac{n-1}{2^{k-1}}l\leq x<\frac{n}{2^{k-1}}l,\\
0,&\text{otherwise},
\end{array}
\right.
\end{equation}
\begin{equation}\label{2.10}
\psi_{n,m}^{''}(x)=\left\{\begin{array}{ll}
\frac{2^{\frac{5k}{2}}}{l^2}\sqrt{\frac{2m+1}{2l}}P_m^{''}\left(\frac{2^k}{l}x-2n+1\right),& \frac{n-1}{2^{k-1}}l\leq x<\frac{n}{2^{k-1}}l,\\
0,&\text{otherwise}.
\end{array}
\right.
\end{equation}

An arbitrary function $f(x)\in L^2[0,l)$ may be approximated using Legendre wavelets as
\begin{equation*}
f(x)\simeq \Psi_{k,M-1}(f) =\sum_{n=1}^{2^{k-1}}\sum_{m=0}^{M-1}a_{n,m}\psi_{n,m}(x)=A^T\psi(x),
\end{equation*}
where
\begin{equation}\label{3.0}
a_{n,m}=\langle f(x),\psi_{n,m}(x)\rangle=\int_0^l f(x)\psi_{n,m}(x)dx,
\end{equation}
\begin{equation*}
A=[a_{1,0},a_{1,1},\ldots,a_{1,M-1},\ldots ,a_{2^{k-1},0},a_{2^{k-1},1},\ldots ,a_{2^{k-1},M-1}]^T,
\end{equation*}
and
\begin{equation*}
\psi(x)=[\psi_{1,0}(x),\psi_{1,1}(x),\ldots ,\psi_{1,M-1}(x),\ldots ,\psi_{2^{k-1},0}(x),\psi_{2^{k-1},1}(x),\ldots ,\psi_{2^{k-1},M-1}(x)]^T.
\end{equation*}
%%%%%%%%%%%%%%%%%%%%%%%%%%%%%%%%%%%%%%%%%%%%%%%%%%%%%%%%%%%%%%%%%%%%%%%%
\section{Numerical method}\label{sec:3}
In order to introduce the numerical method for solving equation (\ref{1.1}), we assume that its solution can be approximated by Legendre wavelets up to $k$-th level of resolution:
\begin{equation}\label{3.1}
y(x)\simeq \sum_{n=1}^{2^{k-1}}\sum_{m=0}^{M-1}a_{n,m}\psi_{n,m}(x),
\end{equation}
where $M\geq 3$. Then, for all $x$ such that $px-\tau \geq 0$, we can write
\begin{equation}\label{3.2}
y(px-\tau)\simeq \sum_{n=1}^{2^{k-1}}\sum_{m=0}^{M-1}a_{n,m}\psi_{n,m}(px-\tau).
\end{equation}
Therefore, if $px-\tau > 0$, by substituting equations (\ref{3.1}) and (\ref{3.2}) in (\ref{1.1}), the following residual is obtained:
\begin{equation*}
\begin{split}
R_1(x)=\sum_{n=1}^{2^{k-1}}\sum_{m=0}^{M-1}a_{n,m}{^cD^\alpha} \psi_{n,m}(x)-f\left(x,\sum_{n=1}^{2^{k-1}}\sum_{m=0}^{M-1}a_{n,m}\psi_{n,m}(x),\right.\\
 \sum_{n=1}^{2^{k-1}}\sum_{m=0}^{M-1}a_{n,m} \psi_{n,m}'(x),\sum_{n=1}^{2^{k-1}}\sum_{m=0}^{M-1}a_{n,m}\psi_{n,m}(px-\tau),\\
\quad \quad\left.\sum_{n=1}^{2^{k-1}}\sum_{m=0}^{M-1}a_{n,m} \psi'_{n,m}(px-\tau),\int_{px-\tau}^xg(x,s,\sum_{n=1}^{2^{k-1}}\sum_{m=0}^{M-1}a_{n,m}\psi_{n,m}(s))ds\right).\\
\end{split}
\end{equation*}
On the other hand, if $px-\tau \leq 0$,  we define the residual as follows:
\begin{equation*}
\begin{split}
R_2(x)=\sum_{n=1}^{2^{k-1}}\sum_{m=0}^{M-1}a_{n,m}{^cD^\alpha} \psi_{n,m}(x)-f\left(x,\sum_{n=1}^{2^{k-1}}\sum_{m=0}^{M-1}a_{n,m}\psi_{n,m}(x),\right.\\
 \sum_{n=1}^{2^{k-1}}\sum_{m=0}^{M-1}a_{n,m} \psi_{n,m}'(x),y_0(px-\tau),y'_0(px-\tau),\\
  \left.\int_{px-\tau}^0g(x,s,y_0(s))ds+\int_{0}^xg(x,s,\sum_{n=1}^{2^{k-1}}\sum_{m=0}^{M-1}a_{n,m}\psi_{n,m}(s))ds\right)
\end{split}
\end{equation*}
The functions $\psi_{n,m}'(\cdot)$ are obtained using \eqref{2.9} and are substituted in $R_1$ and $R_2$. In a particular case, if $\alpha=2$, we have ${^cD^2} \psi_{n,m}(x)=\psi_{n,m}^{''}(x)$ which are given by \eqref{2.10}. Suppose that $1<\alpha<2$, then using \eqref{2.11}, we have
\begin{equation}\label{3.20}
{^cD^\alpha} \psi_{n,m}(x)=\frac{1}{\Gamma(2-\alpha)}\int_{0}^x(x-s)^{1-\alpha}\psi_{n,m}^{''}(s)ds.
\end{equation}
In order to compute the integral part of \eqref{3.20}, first the interval $[0,x]$ is transformed into $[-1,1]$ using the following variable transformation
\begin{equation*}
t=2\left(\frac{s}{x}\right)-1,\quad dt=\frac{2}{x}ds.
\end{equation*}
Therefore, \eqref{3.20} can be rewritten as
\begin{equation}\label{3.21}
{^cD^\alpha} \psi_{n,m}(x)=\frac{1}{\Gamma(2-\alpha)}{\left(\frac{x}{2}\right)}^{2-\alpha}\int_{-1}^{1}(1-t)^{1-\alpha}\psi_{n,m}^{''}(\frac{x}{2}(t+1))dt.
\end{equation}
Since $1-\alpha>-1$, we introduce $\lambda=1-\alpha$ and $\nu=0$ as a special case of Jacobi polynomials and use the Gauss--Jacobi quadrature to compute the integral part of \eqref{3.21}. We set $N:=\lceil \frac{M-2}{2}\rceil$, and get
\begin{equation}\label{3.22}
{^cD^\alpha} \psi_{n,m}(x)=\frac{1}{\Gamma(2-\alpha)}{\left(\frac{x}{2}\right)}^{2-\alpha}\sum_{i=1}^Nw_i\psi_{n,m}^{''}(\frac{x}{2}(t_i+1)),
\end{equation}
where $t_i$, $i=1,\ldots,N$, are the zeros of $J_N^{(1-\alpha,0)}$ and the weights $w_i$, $i=1,\ldots,N$, are given using \eqref{2.8}, as follows:
\begin{equation*}
w_i=\frac{2^{2-\alpha}}{(\frac{d}{dx}J_N^{(1-\alpha,0)}(t_i))^2(1-t_i^2)}.
\end{equation*}
For computing the nodes and weights of the fractional Gauss--Jacobi quadrature numerical algorithms can be used such as Golub-Welsch algorithm introduced by Pang et al. \cite{Pang}.
\begin{remark}\label{rem3.1}
By setting $N:=\lceil \frac{M-2}{2}\rceil$, according to the error of the Gauss--Jacobi quadrature formula, \eqref{3.22} is exact for computing ${^cD^\alpha} \psi_{n,m}(x)$, $n=1$, $m=0,\ldots,M-1$.
\end{remark}
\begin{remark}\label{rem3.2}
In the case of complicated functions $g(\cdot,\cdot,\cdot)$ in $R_1(\cdot)$ and $R_2(\cdot)$, first we transform the interval of integration into the interval $[-1,1]$ and then use the Gauss--Legendre quadrature formula.
\end{remark}

We use the following shifted Chebyshev points as the collocation points
\begin{equation*}
x_{i,j}=\frac{l}{2^k}\left(x_j+2i-1\right),\quad i=1,2,\ldots, 2^{k-1},\quad j=2,3,\ldots , M-1,
\end{equation*}
where
\[
x_j=\cos\left(\frac{(2j+1)\pi}{2M}\right),\quad    j=0,1,\dots,M-1.
\]
Then, set the  residual, $R_l(x)$, $l=1,2$, equal to zero at these points. So that, for  $i=1,2,\ldots, 2^{k-1}$ and $j=2,3,\ldots , M-1$, if $px_{i,j}-\tau > 0$, we set
\begin{equation}\label{3.3}
R_1(x_{i,j})=0,
\end{equation}
otherwise, if $px_{i,j}-\tau \leq 0$, we set
\begin{equation}\label{3.30}
R_2(x_{i,j})=0,
\end{equation}
which \eqref{3.3} and \eqref{3.30} yield a system of $2^{k-1}(M-2)$ nonlinear equations in terms of the $2^{k-1}M$ unknown parameters $a_{n,m}$. To find all the values $a_{n,m}$, it is necessary to complete the system in (\ref{3.3}) with additional conditions in order to obtain a continuous solution. To do this, let us set
\begin{equation}\label{3.4}
\begin{split}
\lim\limits_{x\rightarrow \left(\frac{i-1}{2^{k-1}}l\right)^-}\left(\sum_{n=1}^{2^{k-1}}\sum_{m=0}^{M-1}a_{n,m}\psi_{n,m}(x)\right)-\lim\limits_{x\rightarrow \left(\frac{i-1}{2^{k-1}}l\right)^+}\left(\sum_{n=1}^{2^{k-1}}\sum_{m=0}^{M-1}a_{n,m}\psi_{n,m}(x)\right)=0,\\
i=2,3,\ldots,2^{k-1},~~~~~~~~~~~~
\end{split}
\end{equation}
\begin{equation}\label{3.5}
\begin{split}
\lim\limits_{x\rightarrow \left(\frac{i-1}{2^{k-1}}l\right)^-}\left(\sum_{n=1}^{2^{k-1}}\sum_{m=0}^{M-1}a_{n,m}\psi_{n,m}'(x)\right)-\lim\limits_{x\rightarrow \left(\frac{i-1}{2^{k-1}}l\right)^+}\left(\sum_{n=1}^{2^{k-1}}\sum_{m=0}^{M-1}a_{n,m}\psi_{n,m}'(x)\right)=0,\\
i=2,3,\ldots,2^{k-1}.~~~~~~~~~~~~
\end{split}
\end{equation}
Now, two more equations are needed. These equations can be obtained from the initial conditions (\ref{1.2}) or boundary conditions (\ref{1.3}), respectively, as
\begin{equation}\label{3.6}
\left\{
\begin{split}
&\sum_{n=1}^{2^{k-1}}\sum_{m=0}^{M-1}a_{n,m}\psi_{n,m}(0)=y_0(0),\\
&\lim\limits_{x\rightarrow 0^+}\left(\sum_{n=1}^{2^{k-1}}\sum_{m=0}^{M-1}a_{n,m}\psi_{n,m}'(x)\right)=y'_0,\\
\end{split}
\right.
\end{equation}
\begin{equation}\label{3.7}
\left\{
\begin{split}
&\sum_{n=1}^{2^{k-1}}\sum_{m=0}^{M-1}a_{n,m}\psi_{n,m}(0)=y_0(0),\\
&\lim\limits_{x\rightarrow l^-}\left(\sum_{n=1}^{2^{k-1}}\sum_{m=0}^{M-1}a_{n,m}\psi_{n,m}(x)\right)=y_1.\\
\end{split}
\right.
\end{equation}
Finally, by considering (\ref{3.3})--(\ref{3.5})--(\ref{3.6}) or (\ref{3.3})--(\ref{3.5}),  (\ref{3.7}) in the case of initial or boundary value problem, respectively, we have a system of $2^{k-1}M$ nonlinear equations which can be solved for computing the unknown parameters $a_{n,m}$.
%%%%%%%%%%%%%%%%%%%%%%%%%%%%%%%%%%%%%%%%%%%%%%%%%%%%%%%%%%%%%%%%%%%
\section{Error bounds}\label{sec:4}
Suppose $\mathbb{P}_{M-1}(J_{k,n})$ denotes the set of all functions whose restriction on each mesh interval
$J_{k,n}=(\frac{n-1}{2^{k-1}}l , \frac{n}{2^{k-1}}l)$, $n=1,2,...,{2^{k-1}}$,
are polynomials of degree at most $M-1$. In this section, we give some estimates for the Legendre wavelets truncated error $f(x)-\sum\limits_{n = 1}^{{2^{k - 1}}} {\sum\limits_{m = 0}^{M - 1} {{a_{n,m}}{\psi _{n,m}}(x)} }$ in terms of the Sobolev norms where $\psi _{n,m}(x)$ and  $a_{n,m}$ are given by \eqref{2.2} and \eqref{3.0}, respectively. The Sobolev norm of integer order $m\geq 0$ in the interval (a,b), is given by
\begin{equation}\label{Norm}
\parallel f \parallel _{ H^{m}(a,b)}=\left({\sum^m _{j=0}\parallel f^{(j)}\parallel^2 _{L^2(a,b)}} \right)^\frac{1}{2},
\end{equation}
where $f^{(j)}$ denotes the derivative of $f$ of order $j$. Here, we recall the following results which have been given in \cite{Canuto}.

\begin{lem} \cite{Canuto} \label{Estimates for the Interpolation Error}
Assume that $f(x)\in H^m(-1,1),$ ($H^m(-1,1)$ is a
Sobolev space) with $m\geq 1$ and $L_{M-1}(f) \in \mathbb{P}_{M-1}$ denotes the truncated Legendre series of $f$.
Then the following estimate holds:
\begin{equation}\label{lem1.1}
\parallel{ f-L_{M-1}(f)} \parallel_{L^2(-1,1)}\leq C(M-1)^{-m}|f|_{H^{m,M-1}(-1,1) },
\end{equation}
and for $1\leq r \leq m$, we have
\begin{equation}\label{lem1.2}
\parallel{ f-L_{M-1}(f)} \parallel_{H^r(-1,1)}\leq C(M-1)^{2r-\frac{1}{2}-m}|f|_{H^{m,M-1}(-1,1) },
\end{equation}
where
\begin{equation}\label{Lemma1.2}
\mid f \mid _{ H^{m,M-1}(-1,1)}=\left({\sum^m _{j=min\{m,M\}}\parallel f^{(j)}\parallel^2 _{L^2(-1,1)}} \right)^\frac{1}{2},
\end{equation}
and $C$ is a positive constant independent of the function $f$ and integer $M$.
\end{lem}
\begin{lem} \label{lem2-1}Let $f_{n,k} : J_{k,n} \rightarrow \mathbb{R}$, $n = 1, 2,... ,2^{k-1}$, be functions
in $H^m(J_{k,n})$. Consider the function $\bar{f}_{n,k}: (-1,1) \rightarrow  \mathbb{R}$ defined by
$(\bar{f}_{n,k})(x)=f_{n,k}(\frac{l}{2^{k}}(x+2n-1))$ for all $x \in (-1,1)$, then for $0\leq j\leq m$, we have
\begin{equation}\label{lemma2}
\parallel (\bar{f}_{n,k})^{(j)} \parallel _{ L^2(-1,1)}=({\frac{2^{k}}{l}})^{\frac{1}{2}-j}\parallel f^{(j)}_{n,k}\parallel _{L^2(J_{k,n})}.
\end{equation}
\end{lem}
\textbf{Proof.} For $0\leq j \leq m$, using the definition of the $L^2$- norm and the change of variable by setting $x'=\frac{l}{2^{k-1}}(x+2n-1)$, we have
\begin{equation*}
\begin{split}
\parallel {{{{\bar {f}_{n,k}}}^{(j)}}} \parallel_{{L^2}( - 1,1)}^2 &= \int_{ - 1}^1  |{\bar{f}_{n,k}}^{(j)} (x)|^2dx \\
&=\int_{ - 1}^1 |{f_{n,k}}^{(j)} (\frac{{l}}{2^{k}}(x+2n-1)|^2dx \\
&= \int_{\frac{{n - 1}}{2^{k-1}}l}^{\frac{n}{2^{k-1}}l} {(\frac{{2^{k}}}{{{l}}}} {)^{ - 2j}}{| {f_{n,k}^{(j)}(x')} |^2}(\frac{{2^{k}}}{{{l}}})dx'\\
&= {(\frac{{2^{k}}}{{{l}}})^{1 - 2j}}\parallel {f_{n,k}^{(j)}}\parallel_{{L^2}(J_{k,n})}^2.
\end{split}
\end{equation*}

Now we introduce the following semi-norm of $f \in H^m(0,l)$:
\begin{equation}\label{semi1}
\mid f \mid _{ H^{r,m,M-1,k}(0,l)}=\left({\sum^m _{j=min\{m,M\}}
(2^{k})^{2r-2j} \parallel f^{(j)}
\parallel ^2_{L^2(0,l)}} \right)^\frac{1}{2},
\end{equation}
for $0\leq r \leq m$, $M\geq 1$ and $k\geq 1$. Note that for $M \geq m$, we have
\begin{equation}\label{semi2}
\mid f \mid _{ H^{r,m,M-1,k}(0,l)}=
(2^{k})^{r-m} \parallel f^{(m)}
\parallel _{L^2(0,l)}.
\end{equation}
%%%%%%%%%%%%%%%%%%%%%%%%%%%%%%%%%%%%%%%%%%%%%%%%%%%%%%%%%%%%%%%%%%%%%%%%%%%%%%%%%%%%%%%%%%%%%%%%%%%%%%%%%%%%%%%%%%%%%%%
\begin{theorem} Suppose that $f \in H^m(0,l)$ with $m \geq 1$ and $\linebreak{\Psi_{k,M-1}(f)=\sum\limits_{n = 1}^{{2^{k - 1}}} {\sum\limits_{m' = 0}^{M - 1} {{a_{n,m'}}{\psi _{n,m'}}(x)} }}$ is the best approximation of $f$ based on the Legendre wavelets. Then,
\begin{equation}\label{Theorem}
\parallel{f - \Psi_{k,M-1}(f)} \parallel_{{L^2}(0,l)}\leq C(M-1)^{-m}|f|_{H^{0,m,M-1,k}(0,l)},
\end{equation}
and for $1\leq r \leq m$,
\begin{equation}\label{Theoremm}
\parallel{f - \Psi_{k,M-1}(f)} \parallel_{{H^r}(0,l)}\leq C(M-1)^{2r-\frac{1}{2}-m}|f|_{H^{r,m,M-1,k}(0,l)},
\end{equation}
where $C$ denotes a positive constant that is independent of $M$ and $k$ but depends on the length $l$.
\end{theorem}
%%%%%%%%%%%%%%%%%%%%%%%%%%%%%%%%%%%%%%%%%%%%%%%%%%%%%%%%%%%%%%%%%%%%%%%%%%%%%%%%%%%%%%%%%%%%%%%%%%%%%%%%%%%%%%%%%%%%%%
\textbf{Proof.} Consider the function $f_{n,k} : J_{k,n}\rightarrow \mathbb{R}$
such that $f_{n,k}(x)=f(x)$, for all $x \in J_{k,n}$. Then from equation (\ref{Norm}) and Lemma \ref{lem2-1} for $r \geq 0$ we have \\
\begin{equation*}
\begin{split}
\left\| {{f} - {\Psi _{k,M - 1}}({f})} \right\|_{H^r(0,l )}^2 &= \sum\limits_{n = 1}^{{2^{k - 1}}} {\left\| {{f_{n,k}} - \sum\limits_{m'=0}^{M - 1} {{a_{n,m'}}{\psi _{n,m'}}(x)} } \right\|} _{H^r({J_{k,n}})}^2\\ 
&= c_1\sum\limits_{n = 1}^{2^{k - 1}} \sum\limits_{j = 0}^r {(2^{k})}^{2j - 1}{\left\| {\bar f}_{n,k}^{(j)} -{\left( {\sum\limits_{m' = 0}^{M - 1} \frac{\langle{\bar f}_{n,k},P_{m'}\rangle} {\langle P_{m'},P_{m'}\rangle}{P_{m'}} } \right)}^{(j)} \right\|}_{L^2( - 1,1)}^2,
\end{split}
\end{equation*}
Setting $r=0$, we obtain
\begin{equation*}
\begin{split}
\left\| {{f} - {\Psi _{k,M - 1}}({f})} \right\|_{L^2(0,l )}^2 &= {c_1} \times \frac{1}{{{2^{k}}}}\sum\limits_{n = 1}^{{2^{k - 1}}} {\left\| {{{\bar f}_{n,k}} - \sum\limits_{m' = 0}^{M - 1} {\frac{\langle{\bar f}_{n,k},P_{m'}\rangle} {\langle P_{m'},P_{m'}\rangle}{P_{m'}}} } \right\|} _{L^2( - 1,1)}^2\\
& \le {c_2} \times \frac{1}{{{2^{k}}}}{(M - 1)^{ - 2m}}\sum\limits_{n = 1}^{{2^{k - 1}}} {\sum\limits_{j = \min \{ m,M\} }^m {\left\| {\bar f_{n,k}^{(j)}} \right\| _{L^2( - 1,1)}^2} }\\
& \le {c_3}{(M - 1)^{ - 2m}}\sum\limits_{j = \min \{ m,M\} }^m {{{({2^{k}})}^{ - 2j}}\sum\limits_{n = 1}^{{2^{k - 1}}} {\left\| {f_{n,k}^{(j)}} \right\|^2 _{L^2({J_{k,n}})}} }\\
&= {c_3}{(M - 1)^{ - 2m}}\sum\limits_{j = \min \{ m,M\} }^m {{{({2^{k}})}^{ - 2j}}\left\| {f_{n,k}^{(j)}} \right\|_{L^2(0,l )}^2} ,
\end{split}
\end{equation*}
here, we used equation \eqref{lem1.1} for the second inequality and Lemma \ref{lem2-1} for the third inequality. This yields equation (\ref{Theorem}). Furthermore, for $1\leq r \leq m$ and $k\geq 1$, the following inequality is observed:
\begin{equation*}
\begin{split}
\left\| {{f} - {\Psi _{k,M - 1}}({f})} \right\|_{H^r(0,l)}^2 &= \sum\limits_{n = 1}^{{2^{k - 1}}} {\left\| {{f_{n,k}} - \sum\limits_{m' = 0}^{M - 1} {{a_{n,m'}}} {\psi _{n,m'}}(x)} \right\|} _{H^r({J_{k,n}})}^2\\
&= {c_1}\sum\limits_{n = 1}^{{2^{k - 1}}} {\sum\limits_{j = 0}^r {{{({2^{k}})}^{2j - 1}}\left\| {\bar f_{n,k}^{(j)} - {{\left( {\sum\limits_{m' = 0}^{M - 1} {\frac{\langle{\bar f}_{n,k},P_{m'}\rangle} {\langle P_{m'},P_{m'}\rangle}{P_{m'}}} } \right)}^{(j)}}} \right\|} } _{L^2( - 1,1)}^2\\
& \le {c_2}{(M - 1)^{4r-1 - 2m}}{({2^{k}})^{2r - 1}}\sum\limits_{n = 1}^{{2^{k - 1}}} {\sum\limits_{j = \min \{ m,M\} }^m {\left\| {\bar f_{n,k}^{(j)}} \right\| _{L^2( - 1,1)}^2} }\\
&\le {c_3}{(M - 1)^{4r -1- 2m}}{({2^{k}})^{2r - 1}}\sum\limits_{j = \min \{ m,M\} }^m {{{({2^{k}})}^{1 - 2j}}\sum\limits_{n = 1}^{{2^{k - 1}}} {\left\| {f_{n,k}^{(j)}} \right\| _{L^2({J_{k,n}})}^2} }\\
& = {c_3}{(M - 1)^{4r -1- 2m}}\sum\limits_{j = \min \{ m,M\} }^m {{{({2^{k}})}^{2r - 2j}}\left\| {f_{n,k}^{(j)}} \right\|_{L^2(0,l )}^2} ,
\end{split}
\end{equation*}
where for the third inequality, equation (\ref{lem1.2}) is used. Therefore, we have proved (\ref{Theoremm}).
\begin{theorem}\label{th2} 
Suppose that $y(x)\in H^m(0,l)$ is the exact solution of equation \eqref{1.1}, $\Psi_{k,M-1}(y)=\sum\limits_{n=1}^{2^{k-1}}\sum\limits_{m'=0}^{M-1}a_{n,m'}\psi_{n,m'}(x)=A^T\psi(x)$ is its best approximation using the Legendre wavelets ($M\geq 3$) and $y_{k,M-1}(x)=\sum\limits_{n=1}^{2^{k-1}}\sum\limits_{m'=0}^{M-1}\bar{a}_{n,m'}\psi_{n,m'}(x)=\bar{A}^T\psi(x)$ is the approximate solution obtained by the proposed method in this paper. Then, we have
\begin{equation}\label{8.1}
|| y - y _{k,M - 1}||_{L^2(0,l )}\leq C(M-1)^{-m}|y|_{H^{0,m,M-1,k}(0,l)}+M^{\frac{1}{2}}2^{\frac{k-1}{2}}||A-\bar{A}||_2,
\end{equation}
where the norm in $||A-\bar{A}||_2$  is the usual Euclidean norm for vectors.
\end{theorem}
\textbf{Proof.} Using the triangle property of $L^2$-norm and \eqref{Theorem}, we obtain
\begin{equation}\label{8.2}
\begin{split}
|| y - y _{k,M - 1}||_{L^2(0,l )}&\leq || y - \Psi_{k,M-1}(y)||_{L^2(0,l )}+|| \Psi_{k,M-1}(y)-y _{k,M - 1}||_{L^2(0,l )}\\
& \leq C(M-1)^{-m}|y|_{H^{0,m,M-1,k}(0,l)}+|| \Psi_{k,M-1}(y)-y _{k,M - 1}||_{L^2(0,l )}.
\end{split}
\end{equation}
We have 
\begin{equation}\label{8.3}
\begin{split}
|| \Psi_{k,M-1}(y)-y _{k,M - 1}||_{L^2(0,l )}^2&=\int_0^l\left( \sum_{n=1}^{2^{k-1}}\sum_{m'=0}^{M-1}(a_{n,m'}-\bar{a}_{n,m'})\psi_{n,m'}(x)\right)^2dx\\
&\leq \int_0^l \left(\sum_{n=1}^{2^{k-1}}\sum_{m'=0}^{M-1}|a_{n,m'}-\bar{a}_{n,m'}|^2\right)\left(\sum_{n=1}^{2^{k-1}}\sum_{m'=0}^{M-1}|\psi_{n,m'}(x)|^2\right)dx\\
& = \left(\sum_{n=1}^{2^{k-1}}\sum_{m'=0}^{M-1}|a_{n,m'}-\bar{a}_{n,m'}|^2\right)\left(\sum_{n=1}^{2^{k-1}}\sum_{m'=0}^{M-1}\int_0^l |\psi_{n,m'}(x)|^2dx\right)\\
& =||A-\bar{A}||_2^2 M2^{k-1}.\\
\end{split}
\end{equation}
Taking \eqref{8.2} and \eqref{8.3} into consideration, we obtain \eqref{8.1}.
\begin{remark}\label{rem1}
Theorem \ref{th2} suggests that if the values of $M$ and $k$ are large simultaneously, the numerical solution obtained by the proposed method in this paper would be divergent, specially when the exact solution $y(x)$ belongs to the space $H^m(0,l)$ with small $m$. 
\end{remark}
%%%%%%%%%%%%%%%%%%%%%%%%%%%%%%%%%%%%%%%%%%%%%%%%%%%%%%%%%%%%%%%%%%%%%%%%%%%%%%%%%%%%%
\section{Numerical examples}\label{sec:5}
\begin{figure}[!ht]
\centering
\includegraphics[scale=0.7]{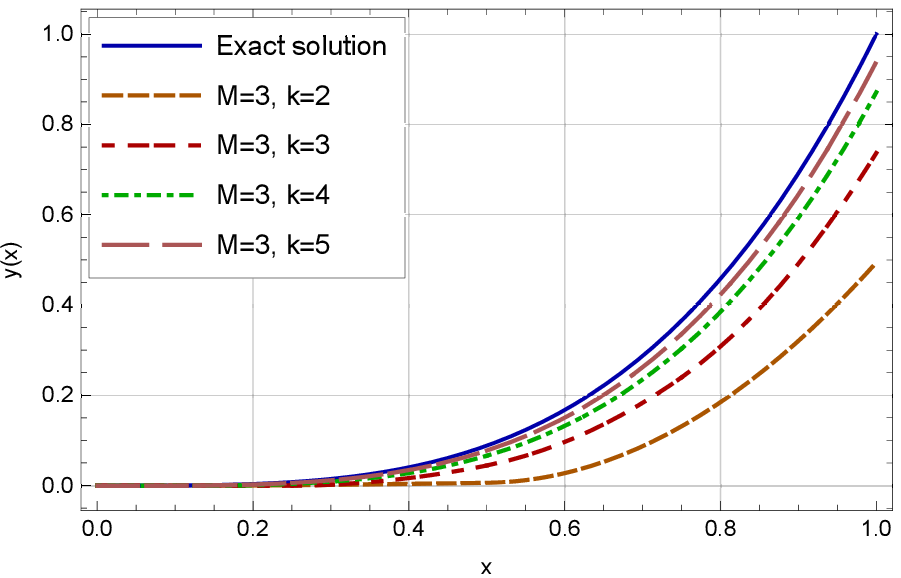}
\includegraphics[scale=0.7]{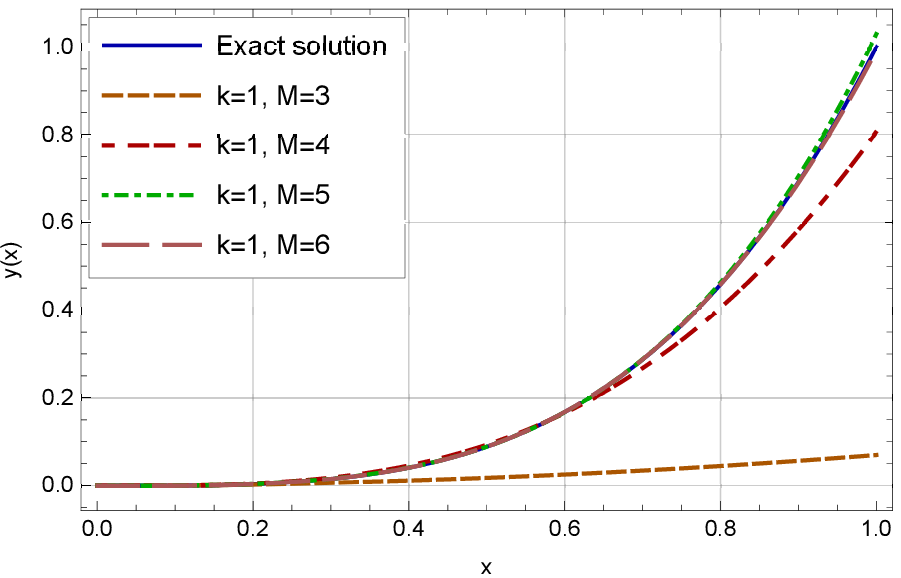}
\caption{(Example \ref{ex2}.) The comparison between the exact solution and approximate solutions with $M=3$ and $k=2,3,4,5$ (left), and with $k=1$ and $M=3,4,5,6$ (right).}\label{fig:4}
\end{figure}

In this section, we consider some examples and solve them using the present method. In our implementation, the method was carried out using \textsf{Mathematica 11.3}. The error $||{y - y_{k,M-1}}||$ is computed in the $L^2$-norm, where $y$ is the exact solution and $y_{k,M-1}$ is the numerical solution obtained by the present method. 
%%%%%%%%%%%%%%%%%%%%%%%%%%%%%%%%%%%%%%%%%%%%%%%%%%%%%%%%%%%%%%%%
\begin{ex}\label{ex1} As the first test example, we consider the following nonlinear fractional delay integro-differential equation
\begin{equation*}
^cD^{1.5} y(x)=y(x-1)y^{'}(x)-\int_{x-1}^x y(s)ds-2x^3+5x^2-3x+\frac{4}{\sqrt{\pi}}x^{\frac{1}{2}}+\frac{1}{3},~\quad 0\leq x\leq 2, 
\end{equation*}
with conditions
\begin{equation*}
y(x)=x^2,\quad \text{if $x\leq 0$},\quad y(2)=4.
\end{equation*}
\end{ex}
The exact solution is $y(x)=x^2$.  This very simple example, that can be solved manually, illustrates how the proposed method can be applied to a boundary value problem. We employ the method with $k=1$, $M=3$, and $l=2$. So, we have
\begin{equation}\label{8.1}
y(x)=a_{1,0}\psi_{1,0}(x)+a_{1,1}\psi_{1,1}(x)+a_{1,2}\psi_{1,2}(x),
\end{equation}
where
\begin{equation*}
\psi_{1,0}(x)=\frac{1}{\sqrt{2}},\quad \psi_{1,1}(x)=\sqrt{\frac{3}{2}} (x-1),\quad \psi_{1,2}(x)=\frac{1}{2} \sqrt{\frac{5}{2}} \left(3 x^2-6 x+2\right).
\end{equation*}
The first and second derivatives of the functions $\psi_{1,m}$, $i=0,1,2$, are given by
\begin{equation*}
\psi_{1,0}^{'}(x)=0,\quad \psi_{1,1}^{'}(x)=\sqrt{\frac{3}{2}},\quad \psi_{1,2}^{'}(x)=\frac{1}{2} \sqrt{\frac{5}{2}} (6 x-6),
\end{equation*}
\begin{equation*}
\psi_{1,0}^{''}(x)=\psi_{1,1}^{''}(x)=0,\quad \psi_{1,2}^{''}(x)=3 \sqrt{\frac{5}{2}}.
\end{equation*}
Moreover, we have
\begin{equation*}
{^cD^{1.5}}\psi_{1,0}(x)={^cD^{1.5}}\psi_{1,1}(x)=0,\quad {^cD^{1.5}}\psi_{1,2}(x)=\frac{1}{\Gamma(0.5)}\int_0^x(x-s)^{-0.5} \left(3\sqrt{\frac{5}{2}}\right)ds.
\end{equation*}
We set $N=1$, and use the Jacobi polynomial $J_1^{(-0.5,0)}(x)=\frac{1}{4} (3 x-1)$ in oreder to employ the Gauss-Jacobi quadrature formula for computing ${^cD^{1.5}}\psi_{1,2}$. The node and weight of the quadrature rule are easily given by
\begin{equation*}
t_1=\frac{1}{3},\quad w_1=2 \sqrt{2}.
\end{equation*}
Hence, using \eqref{3.22}, we obtain
\begin{equation*}
{^cD^{1.5}}\psi_{1,2}=\frac{6 \sqrt{5}}{\Gamma (0.5 )}\left(\frac{x}{2}\right)^{0.5 }.
\end{equation*}
Finally, using the collocation point $x_{1,2}=1-\frac{\sqrt{3}}{2}$ and the boundary conditions given by \eqref{3.7}, the following system of algebraic equations are obtained
\begin{equation*}
\left\{
\begin{split}
&-24 \sqrt{2} \left(\sqrt{3}-2\right) a_{1,0}-42 \sqrt{6} a_{1,1}+\left(57 \sqrt{30}+288 \sqrt{\frac{5 \left(2-\sqrt{3}\right)}{\pi }}\right) a_{1,2}\\
&\quad\quad\quad\quad\quad\quad\quad\quad\quad\quad\quad\quad\quad\quad\quad\quad=12 \sqrt{3}+92 \sqrt{\frac{4-2 \sqrt{3}}{\pi }}-40\\
&\frac{1}{\sqrt{2}}a_{1,0}-\sqrt{\frac{3}{2}} a_{1,1}+\sqrt{\frac{5}{2}} a_{1,2}=0\\
&\frac{1}{\sqrt{2}}a_{1,0}+\sqrt{\frac{3}{2}} a_{1,1}+\sqrt{\frac{5}{2}} a_{1,2}=4.
\end{split}
\right.
\end{equation*}
By solving this system, we get
\begin{equation*}
a_{1,0}=\frac{4 \sqrt{2}}{3},\quad a_{1,1}= 2 \sqrt{\frac{2}{3}},\quad a_{1,2}=\frac{2 }{3}\sqrt{\frac{2}{5}}.
\end{equation*}
By substituting these values into \eqref{8.1}, we have $y(x)=x^2$ which is the exact solution.
%%%%%%%%%%%%%%%%%%%%%%%%%%%%%%%%%%%%%%%%%%%%%%%%%%%%%%%%%%%%%
%%%%%%%%%%%%%%%%%%%%%%%%%%%%%%%%%%%%%%%%%%%%%%%%%%%%%%%%%%%%%%%%%%%%%%%%%%%%%%%%%%%%%%%%%%%%%%%%%
\begin{ex}\label{ex2} Consider the following linear fractional delay pantograph equation
\begin{equation*}
^cD^{1.9} y(x)=y^{'}(\frac{x}{2})-\frac{7 }{8 \sqrt{2}}x^{5/2}+\frac{105 \sqrt{\pi }}{16 \Gamma \left(\frac{13}{5}\right)} x^{8/5},~\quad 0\leq x\leq 1, 
\end{equation*}
\end{ex}
with the initial conditions $y(0)=y'(0)=0$. The exact solution of this problem is $y(x)=x^{7/2}$. We have solved this problem with different values of $k$ and $M$. The numerical results are shown in  Table \ref{tab:5} and Figure \ref{fig:4}. In Table \ref{tab:5} the error norm is displayed for different values of $k$ and $M$.  In the case $M=3$  when $k$ is increased we see that the error norm tends to $0$ like $2^{(-k)}$ (as $k \rightarrow \infty$), which agrees with the approximation error by a Legendre wavelet when the approximated function belongs to $H^3([0,1])$. For higher values of $M$ the convergence rate is even faster.  In Figure \ref{fig:4}, the exact solution and the approximate solutions given by $M=3$, $k=2,3,4,5$  (left) and $k=1$, $M=3,4,5,6$ (right) are displayed.
%%%%%%%%%%%%%%%%%%%%%%%%%%%%%%%%%%%%%%%%%%%%%%%%%%%%%%%%%%%%%%%%%%%%%%%%%%%%%%%%%%%%%%%%
\begin{table}[!ht]
\centering
\caption{(Example \ref{ex2}.) Numerical results with different values of $M$ and $k$.}\label{tab:5}
\begin{tabular}{llllllllll}
\hline
&\multicolumn{2}{c}{$M=3$} && \multicolumn{2}{c}{$M=4$}&&\multicolumn{2}{c}{$M=5$}\\
\cline{2-3}\cline{5-6}\cline{8-9}
$k$ & $L^2$-Error & Ratio && $L^2$-Error & Ratio &&$L^2$-Error & Ratio\\
\hline
$1$ &$3.23e-01$ & --- &  & $5.53e-02$ & ---   &&$7.07e-03$ & --- \\
$2$ &$1.99e-01$ & $1.62$&& $4.73e-02$ &$1.17$ &&$8.06e-03$ & $0.88$\\
$3$ &$1.06e-01$ & $1.88$&& $2.08e-02$ &$2.27$ &&$1.91e-03$ & $4.22$\\
$4$ &$5.23e-02$ & $2.03$&& $7.04e-03$ &$2.95$ &&$4.01e-04$ & $4.76$\\
$5$ &$2.50e-02$ &$2.09$ && $2.04e-03$ &$3.45$ &&$7.60e-05$ & $5.28$\\
$6$ &$1.18e-02$ &$2.19$ && $5.40e-04$ &$3.78$ &&$1.36e-05$ & $5.59$\\
\hline
\end{tabular}
\end{table}

%%%%%%%%%%%%%%%%%%%%%%%%%%%%%%%%%%%%%%%%%%%%%%%%%%%%%%%%%%%%%%%%
\begin{ex}\label{ex3} Consider the following nonlinear fractional pantograph equation \cite{Saeed}
\begin{equation*}
^cD^\alpha y(x)=\frac{8}{3}y'(\frac{x}{2})y(x)+8x^2y(\frac{x}{2})-\frac{4}{3}-\frac{22}{3}x-7x^2-\frac{5}{3}x^3,~ 0\leq x\leq 1,~ 1<\alpha\leq 2,
\end{equation*}
subject to the boundary conditions $y(0) = y(1) = 1$. The exact solution, when $\alpha = 2$, is $y(x) = 1 + x -x^3$.
\end{ex}

\begin{figure}[!ht]
\centering
\includegraphics[scale=0.7]{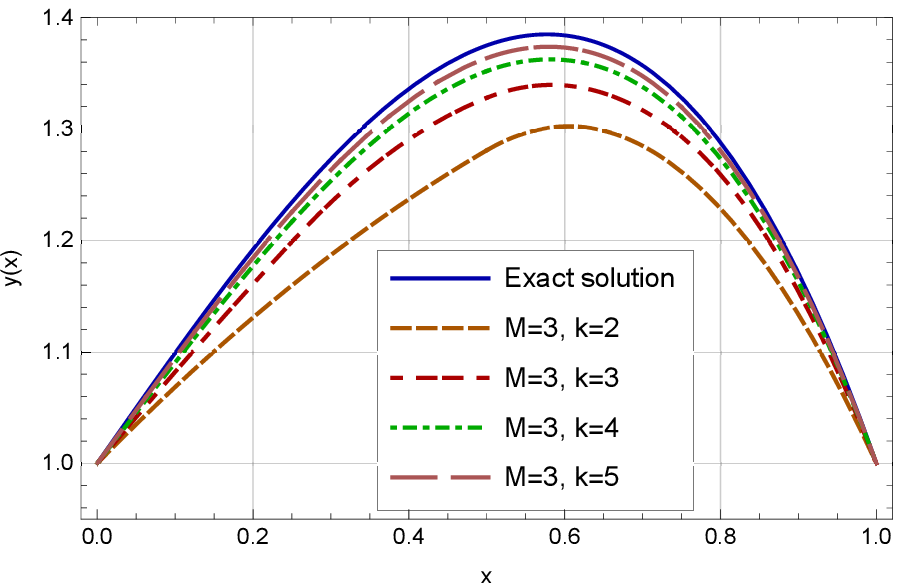}
\includegraphics[scale=0.7]{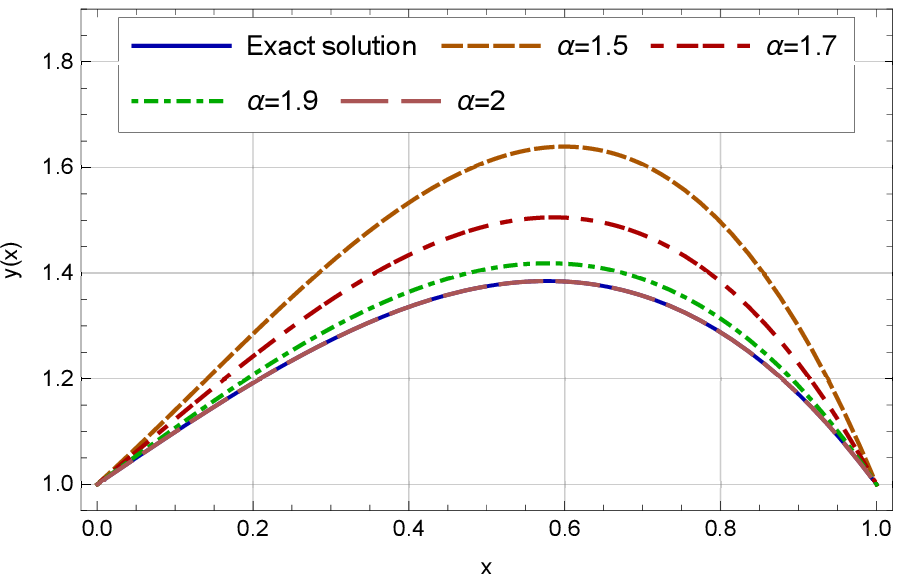}
\caption{(Example \ref{ex3}.) The comparison between the exact solution for $\alpha=2$ and the approximate solutions with $M=3$, $k=2,3,4,5$, $\alpha=2$ (left), and with $M=4$, $k=1$, $\alpha=1.5,1.7,1.9,2$ (right).}\label{fig:2}
\end{figure}

We have solved this problem with different values of $k$ and $M$.  We consider $M=3$, $k=1,2,\ldots,6$, $\alpha=2$ and display the numerical results for the errors in Table \ref{tab:2}.  A plot of the numerical solutions with $M=3$, $k=2,3,4,5$, $\alpha=2$ together with the exact solution is displayed in Figure \ref{fig:2} (left). We can see the convergence of the numerical solutions to the exact solution as the value of $k$ increases. The exact solution for the case $\alpha=2$ and the approximate solutions with $M=4$, $k=1$, and different values of $\alpha$ are seen in Figure \ref{fig:2} (right). As it could be expected, when $\alpha$ is close to $2$, the approximate solution gets close to the exact solution of the case $\alpha=2$. Note that in the case $\alpha=2$, since the exact solution is a third-degree polynomial, our method produces the exact solution with $M\geq 4$. In the case $M=3$, the error tends to zero approximately as $2^{(-k)}$ (as $k \rightarrow \infty$), as it happens in the previous example.
%%%%%%%%%%%%%%%%%%%%%%%%%%%%%%%%%%%%%%%%%%%%%%%%%%%%%%%%%%%%
\begin{table}[!ht]
\centering
\caption{(Example \ref{ex3}.) Numerical results with $M=3$, $\alpha=2$ and different values of $k$. }\label{tab:2}
\begin{tabular}{lll}
\hline
$k$ & $L^2$-Error & Ratio \\
\hline
$1$ &$1.76e-01$ & --- \\
$2$ &$6.89e-02$ & $2.55$\\
$3$ &$3.43e-02$ & $2.01$\\
$4$ &$1.68e-02$ & $2.04$ \\
$5$ &$8.34e-03$ & $2.01$\\
$6$ &$4.15e-03$ & $2.01$\\
\hline
\end{tabular}
\end{table}
%%%%%%%%%%%%%%%%%%%%%%%%%%%%%%%%%%%%%%%%%%%%%%%%%%%%%%%%%%%

\begin{figure}[!ht]
\centering
\includegraphics[scale=0.7]{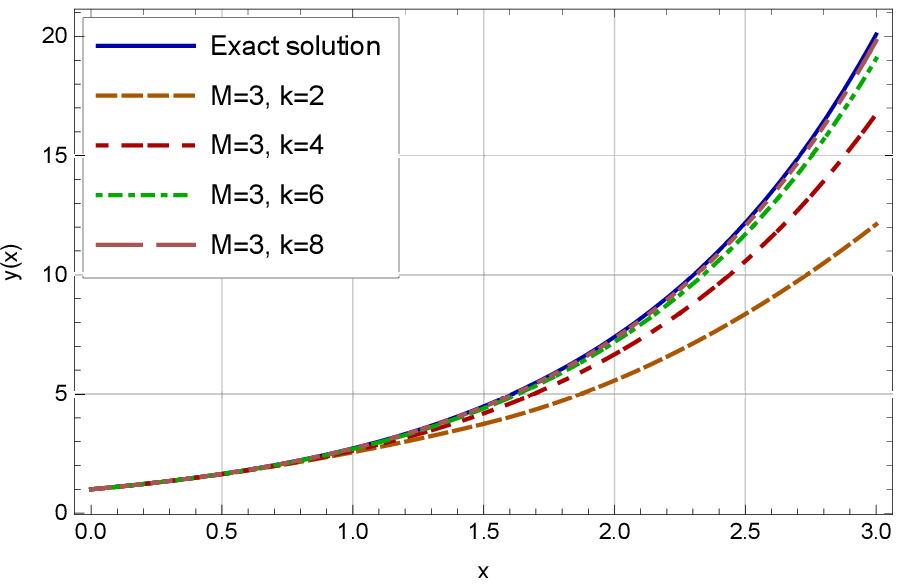}
\includegraphics[scale=0.7]{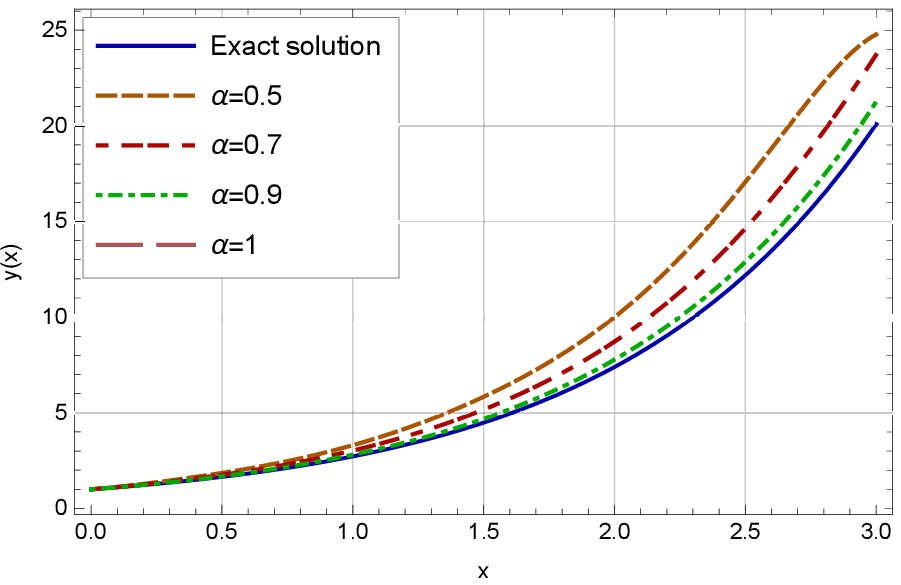}
\caption{(Example \ref{ex4}.) The comparison between the exact solution and approximate solutions with $\alpha=1$, $M=3$, $k=2,4,6,8$ (left) and $M=8$, $k=1$, $\alpha=0.5,0.7,0.9,1$ (right). }\label{fig:3}
\end{figure}

%%%%%%%%%%%%%%%%%%%%%%%%%%%%%%%%%%%%%%%%%%%%%%%%%%%%%%%%%%%%%%%%%%%%%%%%%%%%%%
\begin{ex}\label{ex4} Consider the following fractional Volterra delay integro-differential equation \cite{Saeed}
\begin{equation}\label{5.1}
\begin{split}
&^cD^\alpha y(x)=y(x-1)+\int_{x-1}^x y(s)ds,~ x\geq 0,~ 0<\alpha\leq 1,\\
&y(x)=e^x,\quad\text{if }x\leq 0.
\end{split}
\end{equation}
The exact solution of this problem, when $\alpha = 1$, is $y(x) =e^ x$.
\end{ex}
Considering $\alpha=1$ and substituting $x = 0$ in (\ref{5.1}), we get $y'(0)=1$. By taking the derivative of (\ref{5.1}), we obtain the following
equation
\begin{equation*}
\begin{split}
&^cD^{\alpha+1} y(x)=y'(x-1)+y(x)-y(x-1),~ x\geq 0,~ 0<\alpha\leq 1,\\
&y(x)=e^x,\quad \text{if }x\leq 0.
\end{split}
\end{equation*}
We consider $l=3$ and solve this problem.  A comparison between the absolute error given by our method and the first kind Chebyshev wavelet method introduced in \cite{Saeed} is displayed in Table \ref{tab:3} for $\alpha=1$ and with $k=1$, $M=10, 20$. Furthermore, Table \ref{tab:4} demonstrates the relationship of the errors with the change of $M$ with considering $k=1,2,\ldots,6$ for the classical non-fractional equations corresponding to Example \ref{ex4}. The results of this table suggest that the error tends to zero very fast as $k$ and $M$ tend to infinity. This could be expected since according to Theorem \ref{th2} when a given function is sufficiently smooth the approximation by the Legendre wavelets has spectral convergence. The comparison between the exact solution and approximate solutions obtained in the case $M=3$ and different values of $k$ for $\alpha=1$ is shown in Figure \ref{fig:3} (left). From this figure, it can be seen that the numerical solutions converge to the exact solution as $k$ increases. Finally, the numerical solutions resulted by different values of $\alpha$ and $M=8$, $k=1$ together with the exact solution of the case $\alpha=1$ are plotted in Figure \ref{fig:3} (right). As it could be seen, the approximate solution gets close to the exact solution of the case $\alpha=1$ when $\alpha$ is close to $1$. 
%%%%%%%%%%%%%%%%%%%%%%%%%%%%%%%%%%%%%%%%%%%%%%%%%%%%%%%%%%%%%%%%%%%%
\begin{table}[!ht]
\centering
\caption{(Example \ref{ex4}.) Absolute error at some selected points with $\alpha=1$ and $k=1$.}\label{tab:3}
\begin{tabular}{lllllll}
\hline
 &  \multicolumn{2}{c}{Present method}  & && \multicolumn{2}{c} {Method of \cite{Saeed}}   \\
\cline{2-3}\cline{6-7}
$x$ &   $M=10$ & $M=20$ & && $M=10$ & $M=20$\\
\hline
$0.0$ & $0.00e+00$  & $0.00e+00$&  && $0.00e+00$  & $0.00e+00$\\
$0.6$ & $1.21e-08$   &$3.33e-21$ &&& $1.46e-05$ & $2.13e-13$\\
$1.2$ & $2.45e-08$   &$2.70e-21$ &&& $2.25e-05$  & $3.73e-14$\\
$1.8$ & $1.79e-07$   &$4.93e-21$ &&& $4.21e-05$  & $2.93e-13$\\
$2.4$ & $2.18e-07$  &$3.71e-20$ &&& $7.66e-05$  & $4.23e-13$\\
$3.0$ & $4.35e-05$   &$1.67e-17$ && & $1.45e-04$   & $8.17e-13$\\
\hline
\end{tabular}
\end{table}
%%%%%%%%%%%%%%%%%%%%%%%%%%%%%%%%%%%%%%%%%%%%%%%%%%%%%%%%%%%%%%%%%%%%%%%%%%%%%%%%%%%%%%%%
\begin{table}[!ht]
\centering
\caption{(Example \ref{ex4}.) Numerical results with $\alpha=1$ and different values of $M$ and $k$.}\label{tab:4}
\begin{tabular}{llllllllll}
\hline
&\multicolumn{2}{c}{$M=3$} && \multicolumn{2}{c}{$M=4$}&&\multicolumn{2}{c}{$M=5$}\\
\cline{2-3}\cline{5-6}\cline{8-9}
$k$ & $L^2$-Error & Ratio && $L^2$-Error & Ratio &&$L^2$-Error & Ratio\\
\hline
$1$ &$5.93e+00$ & --- &  & $1.85e+00$ & ---   &&$4.43e-01$ & --- \\
$2$ &$4.63e+00$ & $1.28$&& $1.14e+00$ &$1.62$ &&$1.88e-01$ & $2.36$\\
$3$ &$3.18e+00$ & $1.46$&& $4.86e-01$ &$2.35$ &&$4.23e-02$ & $4.44$\\
$4$ &$1.94e+00$ & $1.64$&& $1.48e-01$ &$3.28$ &&$6.85e-03$ & $6.17$\\
$5$ &$1.10e+00$ & $1.76$&& $4.28e-02$ &$3.46$ &&$9.76e-04$ &$7.02$\\
$6$ &$5.89e-01$ & $1.88$&& $1.11e-02$ &$3.86$ &&$1.29e-04$&$7.57$\\
\hline
\end{tabular}
\end{table}
%%%%%%%%%%%%%%%%%%%%%%%%%%%%%%%%%%%%%%%%%%%%%%%%%%%
%%%%%%%%%%%%%%%%%%%%%%%%%%%%%%%%%%%%%%%%%%%%%%%%%%%%%%%%%%%%%%%%%%%%%%

%%%%%%%%%%%%%%%%%%%%%%%%%%%%%%%%%%%%%%%%%%%%%%
\section{Concluding remarks}\label{sec:6}
In this work, a collocation method based on the Legendre wavelets combined with the Gauss--Jacobi quadrature formula is introduced to solve a class of nonlinear delay-type fractional integro-differential equations. First, an approximation of the unknown function using the Legendre wavelets basis functions is considered. Then, by substituting this approximation in the problem and using the initial condition, we define suitable  residual functions. We use the Gauss--Jacobi quadrature formula in order to compute the fractional derivative of the Legendre wavelets. By collocating the residual at the shifted Gauss--Chebyshev points, a system of nonlinear algebraic equations is obtained. In order to get a continuous solution at the subinterval division nodes, we consider some additional conditions. By using the considered approximation of the solution, the initial or boundary conditions are taken into account. Finally, by solving the resulting system, the unknown coefficients of the solution are computed. Some error bounds are given for the error of the best approximation of a function using the Legendre wavelets functions. To illustrate the efficiency and high accuracy of the suggested method in this paper, five numerical examples have been solved. The numerical examples show the great advantage of using the wavelets method for solving this type of equations. Actually if the solution of the equation is sufficiently smooth by increasing the degree of polynomial ($M$) we can achieve spectral accuracy. But even if the solution has discontinuous derivatives of high order, by choosing moderate value of $M$, fast convergence can be achieved by increasing the level of resolution ($k$).
%%%%%%%%%%%%%%%%%%%%%%%%%%%%%%%%%%%%%%%%%%%%%%%%%%%%%%%%%%%%%%%%%%%%%%%%%%%%%%%%%%%%%%%%%%%%%%%%%%%%%%%%%%%%%%%%%

\section*{Acknowledgements}
P.M. Lima  acknowledges support from Funda\c c\~ao para a Ci\^encia e a Tecnologia (the Portuguese Foundation for Science and Technology) through the grant 
SFRH/BSAB/135130/2017 and UID/Multi/04621/2019.

%%%%%%%%%%%%%%%%%%%%%%%%%%%%%%%%%%%%%%%%%%%%%%%%%%%%%%%%%%%%%%%%%%%%%%%%%%%%%%%%%%%%

%%%%%%%%%%%%%%%%%%%%%%%%%%%%%%%%%%%%%%%%%%%%%%%%%%%%%%%%%%%%%%%%%%%%%%%%%%%%

%%%%%%%%%%%%%%%%%%%%%%%%%%%%%%%%%%%%%%%%%%%%%%%%%%%%%%%%%%%%%%%%%%%%%%%%%%%%%%%%%%%%

\end{document}